\documentclass[12pt,reqno]{amsart}
\usepackage{amssymb,verbatim,enumerate,ifthen}
\usepackage[mathscr]{eucal}
\usepackage[utf8]{inputenc}
\usepackage[T1]{fontenc}
\def\N{\mathbb{N}}
\def\R{\mathbb{R}}

\def\Z{\mathbb{Z}}

\newtheorem{theorem}{Theorem}
\newtheorem*{theorem*}{Theorem}
\def\Thm#1#2{\ifthenelse{\equal{#1}{*}}{\begin{theorem*}#2\end{theorem*}}
             {\begin{theorem}\label{T#1}#2\end{theorem}}}
\newtheorem{Atheorem}{Theorem}

\def\THM#1#2{\begin{Atheorem}\label{T#1}#2\end{Atheorem}}
\def\thm#1{Theorem~\ref{T#1}}

\newtheorem{proposition}[theorem]{Proposition}
\newtheorem*{proposition*}{Proposition}
\def\Prp#1#2{\ifthenelse{\equal{#1}{*}}{\begin{proposition*}#2\end{proposition*}}
             {\begin{proposition}\label{P#1}#2\end{proposition}}}

\newtheorem{corollary}[theorem]{Corollary}
\newtheorem*{corollary*}{Corollary}
\def\Cor#1#2{\ifthenelse{\equal{#1}{*}}{\begin{corollary*}#2\end{corollary*}}
             {\begin{corollary}\label{C#1}#2\end{corollary}}}

\newtheorem{lemma}[theorem]{Lemma}
\newtheorem*{lemma*}{Lemma}
\def\Lem#1#2{\ifthenelse{\equal{#1}{*}}{\begin{lemma*}#2\end{lemma*}}
             {\begin{lemma}\label{L#1}#2\end{lemma}}}
\def\lem#1{Lemma~\ref{L#1}}

\newtheorem{remark}[theorem]{Remark}
\newtheorem*{remark*}{Remark}
\def\Rem#1#2{\ifthenelse{\equal{#1}{*}}{\begin{remark*}\rm #2\end{remark*}}
             {\begin{remark}\label{R#1}\rm #2\end{remark}}}

\newtheorem{example}[theorem]{Example}
\newtheorem*{example*}{Example}
\def\Exa#1#2{\ifthenelse{\equal{#1}{*}}{\begin{example*}\rm #2\end{example*}}
             {\begin{example}\label{Ex#1}\rm #2\end{example}}}

\def\eq#1{{\rm(\ref{E#1})}}
\def\Eq#1#2{\ifthenelse{\equal{#1}{*}}
  {\begin{equation*}\begin{aligned}[]#2\end{aligned}\end{equation*}}
  {\begin{equation}\begin{aligned}[]\label{E#1}#2\end{aligned}\end{equation}}}

\begin{document}
\begin{flushright}
\textit{Submitted to: Acta Sci. Math. (Szeged)} 
\end{flushright}
\vspace{5mm}

\date{\today}

\title[Characterization of generalized quasi-arithmetic means]
{Characterization of generalized quasi-arithmetic means}

\author[J. Matkowski]{Janusz Matkowski}
\address{Faculty of Mathematics, Informatics and Econometrics, 
University of Zielona G\'ora, PL-65-516 Zielona G\'ora, Poland}
\email{j.matkowski@wmie.uz.zgora.pl}

\author[Zs. P\'ales]{Zsolt P\'ales}
\address{Institute of Mathematics, University of Debrecen, 
H-4010 Debrecen, Pf.\ 12, Hungary}
\email{pales@science.unideb.hu}

\dedicatory{Dedicated to the 80th birthday of Professor L\'aszl\'o Leindler}


\thanks{
This research of the second author was realized in the frames of T\'AMOP 4.2.4. A/2-11-1-2012-0001 
”National Excellence Program – Elaborating and operating an inland student and researcher personal
support system”. The project was subsidized by the European Union and co-financed by the European 
Social Fund. This research of the second author was also supported by the Hungarian Scientific 
Research Fund (OTKA) Grant K 111651.}

\begin{abstract} In this paper we characterize generalized quasi-arithmetic means, that is
means of the form $M(x_1,\dots,x_n):=(f_1+\cdots+f_n)^{-1}(f_1(x_1)+\cdots+f_n(x_n))$, where
$f_1,\dots,f_n:I\to\R$ are strictly increasing and continuous functions. Our characterization
involves the Gauss composition of the cyclic mean-type mapping induced by $M$ and a generalized
bisymmetry equation.

\end{abstract}

\maketitle

\section{Introduction}

The notion of quasi-arithmetic mean was introduced in the book of Hardy, Littlewood and P\'olya in 
\cite{HarLitPol34} as a function $A_f:\bigcup_{n=1}^\infty I^n\to I$ defined by
\Eq{Af}{
  A_f(x_1,\dots,x_n):=f^{-1}\Big(\frac{f(x_1)+\cdots+f(x_n)}{n}\Big)\qquad(n\in\N,\,x_1,\dots,x_n\in I),
}
where $I\subseteq\R$ denotes a non-degenerated interval (also in the rest of this paper) and $f:I\to\R$ is a 
continuous strictly monotone function.
The mean $A_f$ is said to be the \emph{quasi-arithmetic mean} generated by $f$. The restriction of $A_f$ to 
$I^n$ will be called the \emph{$n$-variable quasi-arithmetic mean} generated by $f$.

One can easily see that $M=A_f$ is a \emph{mean} in the sense that, for all $n\in\N$ and for all 
$x_1,\dots,x_n\in I$,
\Eq{*}{
  \min(x_1,\dots,x_n)\leq M(x_1,\dots,x_n)\leq \max(x_1,\dots,x_n)
}
holds. Furthermore, $M=A_f$ is a \emph{strict mean} because both inequalities are strict whenever 
$\min(x_1,\dots,x_n)<\max(x_1,\dots,x_n)$.

For the equality problem of quasi-arithmetic means the following result can be established.

\THM{EQ}{Let $f,g:I\to\R$ be continuous and strictly monotone functions. Then the following properties are 
equivalent:
\begin{enumerate}[(i)]
\item For all $n\in\N$ and for all $x_1,\dots,x_n\in I$,
\Eq{*}{
  A_f(x_1,\dots,x_n)=A_g(x_1,\dots,x_n);
}
\item There exists $n\in\N\setminus\{1\}$ such that for all $x_1,\dots,x_n\in I$,
\Eq{*}{
  A_f(x_1,\dots,x_n)=A_g(x_1,\dots,x_n);
}
\item The function $f$ and $g$ are affine transformation of each other, that is, there exist real numbers 
$a,b$ such that $g=af+b$.
\end{enumerate}}

The characterization of quasi-arithmetic means was independently found by Kolmogorov \cite{Kol30}, Nagumo 
\cite{Nag30} and de Finetti \cite{Def31}. The result established by Kolmogorov reads as follows.

\THM{CK}{A function $M:\bigcup_{n=1}^\infty I^n\to I$ is a quasi-arithmetic mean, that is, there exists a 
continuous strictly monotone function $:I\to\R$ such that $M=A_f$ if and only if
\begin{enumerate}[(i)]
 \item for all $n\in\N$, the restriction $M_n:=M|_{I^n}$ is a continuous and symmetric function on $I^n$ 
which is strictly increasing in each of its variables;
 \item for all $n\in\N$, $M_n$ is \emph{reflexive}, that is, $M_n(x,\dots,x)=x$ for all $x\in I$;
 \item $M$ \emph{associative}, that is, for all $n,m\in\N$ and $x_1,\dots,x_n,y_1,\dots,y_m\in I$, we have
\Eq{AS}{
  M_{n+m}(x_1,\dots,x_n,y_1,\dots,y_m)=M_{n+m}(x_1,\dots,x_n,y,\dots,y),
}
where $y=M_m(y_1,\dots,y_m)$.
\end{enumerate}}

The above characterization theorem does not characterize quasi-arithmetic means of fixed number of variables 
because \eq{AS} involves $m$- and $n+m$-variable means. The characterization of $2$-variable 
quasi-arithmetic means was established by Aczél \cite{Acz48a} and this was extended to the $n$-variable case 
by M\"unnich, Maksa and Mokken \cite{MunMakMok99,MunMakMok00}. Their results can be formulated in the 
following way.

\THM{CMMM}{Let $n\geq2$ and let $M:I^n\to I$. Then $M$ is an $n$-variable quasi-arithmetic mean, that is,
$M=A_f|_{I^n}$ for some continuous strictly monotone function $f:I\to\R$ if and only if
\begin{enumerate}[(i)]
 \item $M$ is a continuous and symmetric function on $I^n$ 
which is strictly increasing in each of its variables;
 \item $M$ is reflexive;
 \item $M$ \emph{bisymmetric}, that is, for all $x_{i,j}\in I$ ($i,j\in\{1,\dots,n\}$), we have
\Eq{BS}{
  M(M(x_{1,1},&\dots,x_{1,n}),\dots,M(x_{n,1},\dots,x_{n,n}))\\
  &=M(M(x_{1,1},\dots,x_{n,1}),\dots,M(x_{1,n},\dots,x_{n,n})).
}
\end{enumerate}}

It turns out that also weighted (and therefore, in general, non-symmetric) quasi-arithmetic means can be 
characterized by the bisymmetry equation \eq{BS}.

Quasi-arithmetic means can be generalized in several ways. In 1963 Bajraktarevi\'c \cite{Baj63} introduced 
the notion of quasi-arithmetic means weighted by a weight function. Their equality problem was solved by 
Acz\'el and Dar\'oczy in \cite{AczDar63c}. The characterization theorem of Bajraktarevi\'c means was found
by P\'ales in \cite{Pal87d}. Anoter (symmetric) generalization, the notion of deviation mean, was invented by 
Dar\'oczy in \cite{Dar71b,Dar72b}. The characterization of the Dar\'oczy means was then established by 
P\'ales in \cite{Pal82a}. Both of these characterization theorems use system of functional inequalities 
instead of functional equations like associativity or bisymmetry.

In this paper we consider a recent generalization of quasi-arithmetic means which was introduced by Matkowski 
\cite{Mat10b} in 2010. Given a system $f_1,\dots,f_n:I\to\R$ of continuous strictly increasing functions,
the \emph{generalized $n$-variable quasi-arithmetic mean} $A_{f_1,\dots,f_n}:I^n\to I$ is defined by
\Eq{*}{
  A_{f_1,\dots,f_n}(x_1,\dots,x_n):=(f_1+\cdots+f_n)^{-1}(f_1(x_1)+\cdots+f_n(x_n)) 
   \qquad(x_1,\dots,x_n\in I).
}
The functions $f_1,\dots,f_n$ are called the generators of the mean $A_{f_1,\dots,f_n}$. In the particular 
case $f_1=\cdots=f_n=f$, one can see that $A_{f_1,\dots,f_n}$ reduces to the quasi-arithmetic mean $A_f$. 
More generally, if $f_i=\lambda_i f$, where $\lambda_1,\dots,\lambda_n>0$ and $f:I\to\R$ is a continuous 
strictly increasing function, then $A_{f_1,\dots,f_n}$ will be equal to a so-called weighted quasi-arithmetic
mean. One can easily check that generalized $n$-variable quasi-arithmetic means are strict means.

The equality problem of generalized $n$-variable quasi-arithmetic mean was answered by Mat\-kow\-ski in 
\cite{Mat10b} as follows.

\THM{EM}{Let $f_1,\dots,f_n,g_1,\dots,g_n:I\to\R$ be continuous strictly increasing functions. Then the 
following two assertions are equivalent:
\begin{enumerate}[(i)]
\item For all $x_1,\dots,x_n\in I$,
\Eq{*}{
  A_{f_1,\dots,f_n}(x_1,\dots,x_n)=A_{g_1,\dots,g_n}(x_1,\dots,x_n);
}
\item There exist real numbers $a,b_1,\dots,b_n$ such that, for all $i\in\{1,\dots,n\}$,
\Eq{*}{
  g_i=af_i+b_i.
}
\end{enumerate}}

The main problem addressed in this paper is the characterization of generalized $n$-variable quasi-arithmetic 
means. For this purpose, we recall the notion of Gauss composition in the next section with its basic 
properties and, in the last section, we introduce the family of cyclic mean-type mapping attached to a given 
generalized $n$-variable quasi-arithmetic mean and we compute its Gauss composition explicitly. Using this, 
we shall deduce a bisymmetry type identity for generalized $n$-variable quasi-arithmetic means which will 
turn out to be the key characteristic property beyond regularity and reflexivity properties. The key point in 
our proof is the use of the description of the $CM$-solutions of the generalized bisymmetry equation due to 
Maksa \cite{Mak99b}.

\section{Auxiliary notions and results}

\subsection{Gauss iteration and Gauss composition of means}

Given a system $M=(M_1,\dots,M_n):I^n\to I^n$ of $n$-variable means (which is also called an $n$-variable 
mean-type mapping (cf.\ \cite{Mat99b}) and an element $x\in I^n$, the sequence $x_{k}:=M^k(x)$ is called the 
\emph{Gauss iteration of $x\in I^n$ by the mean-type mapping $M$}.

\THM{GI}{Assume that $M=(M_1,\dots,M_n):I^n\to I^n$ is a continuous strict $n$-variable mean-type mapping, 
that is, $M_1,\dots,M_n$ are continuous and strict means. Then there exists a unique continuous strict mean 
$K:I^n\to I$ such that the sequence of iterates $(M^k)$ converges pointwise to the mean-type mapping 
$(K,\dots,K)$. Furthermore, $K$ is the unique mean satisfying the identity 
\Eq{IE}{
  K\circ M= K,
}
which is called the invariance equation for $K$.}

For the proof of this result, the reader should check the papers \cite{DarPal02c}, \cite{Mat99b}.

The mean $K$ constructed in the above theorem will be called the \emph{Gauss composition of the means 
$(M_1,\dots,M_n)$} and will be denoted by $\Gamma(M_1,\dots,M_n)$.

\subsection{Cyclic mean-type mappings} Given $n\ge2$, define the cyclic permutation
$\sigma:\{1,\dots,n\}\to\{1,\dots,n\}$ by 
\Eq{*}{
\sigma(k):=\begin{cases}n & \mbox{if } k=1, \\ k-1 & \mbox{if } k\in\{2,\dots,n\}.\end{cases}
}
Clearly, $\sigma^n=\sigma^0$ which is the identity mapping. 

In the proof of the main result we shall need the following

\Lem{P}{For all $i\in\{1,\dots,n\}$, we have $\sigma^i(i)=n$.} 

\begin{proof} The statement is obvious for $i=1$ by the definition of $\sigma$.
Assume that it holds for $i=k$, where $k\in\{1,\dots,n-1\}$. Then using this
inductive assumption, we get
\Eq{*}{
  \sigma^{k+1}(k+1)=\sigma^k(\sigma(k+1))=\sigma^k(k)=n.
}
Thus, the statement is also true for $i=k+1$.
\end{proof}

For an $n$-variable mean $M:I^n\to I$ and
an index $i\in\Z$, the $i$th cyclically permuted mean $M^{\langle i\rangle}:I^n\to I$ is defined by
\Eq{*}{
  M^{\langle i\rangle}(x_1,\dots,x_n)
    =M(x_{\sigma^i(1)},\dots,x_{\sigma^i(n)})\qquad(x_1,\dots,x_n\in I),
}
The mapping $(M^{\langle 0\rangle},\dots,M^{\langle n-1 \rangle}):I^n\to I^n$ 
is said to be the \emph{cyclic mean-type mapping induced by $M$}.

\Prp{1}{If $M:I^n\to I$ is a continuous strict $n$-variable mean, then the Gauss 
composition $\Gamma(M^{\langle 0\rangle},\dots,M^{\langle n-1 \rangle})$ is a cyclically 
symmetric mean, i.e., for all $i\in\Z$,
\Eq{*}{
\Gamma(M^{\langle 0\rangle},\dots,M^{\langle n-1 \rangle}) 
  = \big(\Gamma(M^{\langle 0\rangle},\dots,M^{\langle n-1 \rangle})\big)^{\langle i\rangle}.
}}

\begin{proof}
Put $K=\Gamma(M^{\langle 0\rangle},\dots,M^{\langle n-1 \rangle})$. By \thm{GI}, $K$ is the
unique $n$ variable mean which solves the functional equation
\Eq{KI}{
 K\circ(M^{\langle 0\rangle},\dots,M^{\langle n-1 \rangle})=K.
}
Therefore, for $i\in\Z$,
\Eq{*}{
  K^{\langle i\rangle}\circ(M^{\langle 0\rangle},\dots,M^{\langle n-1 \rangle})
  =K\circ(M^{\langle i\rangle},\dots,M^{\langle i+n-1 \rangle})
  =\big(K\circ(M^{\langle 0\rangle},\dots,M^{\langle n-1 \rangle})\big)^{\langle i\rangle}
  =K^{\langle i\rangle}.
}
Thus, $K^{\langle i\rangle}$ is also a solution of the invariance equation \eq{KI}. Hence,
by the unique solvability, it follows that $K=K^{\langle i\rangle}$.
\end{proof}

\subsection{$CM$-quasi-sums, $CM$-functions and generalized bisymmetry} An $n$-variable function 
$F:I^n\to I$ is called a $CM$-\emph{quasi-sum} (cf.\ Maksa \cite{Mak99b}) if there exist $CM$ 
(i.e., continuous and strictly increasing) functions $f_1,\dots,f_n:I\to\R$ and 
$f:f_1(I)+\cdots+f_n(I)\to I$ such that
\Eq{*}{
  F(x_1,\dots,x_n)=f(f_1(x_1)+\cdots+f_n(x_n)) \qquad((x_1,\dots,x_n)\in I^n).
}
The functions $f,f_1,\dots,f_n$ are called the generators of the quasi-sums.
A function $F:I^n\to I$ is said to be a $CM$-function if it is continuous and strictly 
increasing in each of its variables.

The following result, which is a particular case of the general theorem of Maksa \cite{Mak99b}, 
will play a crucial role in our approach.

\THM{MG}{Let $n,m>1$. Let $F,F_1,\dots,F_n:I^m\to I$ and $G,G_1,\dots,G_m:I^n\to I$ be $CM$-functions
such that, for all $x_{i,j}\in I$, ($i\in\{1,\dots,n\}$, $j\in\{1,\dots,m\}$),
\Eq{*}{
  F(G_1(x_{1,1},\dots,x_{n,1}),\dots,G_m(x_{1,m},\dots,x_{n,m}))
  =G(F_1(x_{1,1},\dots,x_{1,m}),\dots,F_n(x_{n,1},\dots,x_{n,m})).
}
Then $F,F_1,\dots,F_n$ and $G,G_1,\dots,G_m$ are $CM$-quasi-sums.}

\Lem{GQ}{A $CM$-quasi-sum $M:I^n\to I$ is reflexive if and only if it is a generalized $n$-variable 
quasi-arithmetic mean.}

\begin{proof} Obviously generalized $n$-variable quasi-arithmetic mean are $CM$-quasi-sums.

Assume now that a $CM$-quasi-sum 
\Eq{*}{
M(x_{1},\dots,x_{n}) =f(f_1(x_1)+\cdots+f_n(x_n)) \qquad (( x_{1},\dots,x_{n}) \in I^{n})
}
is a mean. Setting $x_{1}=\cdots=x_{n}=x,$ by the reflexivity of $M,$ we get 
$f((f_{1}+\cdots+f_n)(x))=x$ for all $x\in I$, whence
\Eq{*}{
f=(f_{1}+\cdots+f_n)^{-1}.
}
Thus, $M$ is of the form $A_{f_1,\dots,f_n}$, i.e., $M$ is generalized $n$-variable quasi-arithmetic mean.
\end{proof}

\section{Characterization of generalized quasi-arithmetic means}

\Thm{M1}{Let $f_1,\dots,f_n:I\to\R$ be continuous strictly increasing functions. Then the Gauss composition
of the means $A^{\langle 0\rangle}_{f_1,\dots,f_n},\dots,A^{\langle n-1\rangle}_{f_1,\dots,f_n}$ is the 
quasi-arithmetic mean $A_{f_1+\dots+f_n}$. Furthermore, for all $x_{i,j}\in I$, ($i,j\in\{1,\dots,n\}$),
the generalized bisymmetry equation
\Eq{IEG}{
  A_{f_1+\dots+f_n}&\big(A^{\langle 0\rangle}_{f_1,\dots,f_n}(x_{1,1},\dots,x_{n,1}),
       \dots,A^{\langle n-1\rangle}_{f_1,\dots,f_n}(x_{1,n},\dots,x_{n,n})\big)\\
  &=A_{f_1+\dots+f_n}\big(A^{\langle 0\rangle}_{f_1,\dots,f_n}(x_{1,1},\dots,x_{1,n}),
       \dots,A^{\langle n-1\rangle}_{f_1,\dots,f_n}(x_{n,1},\dots,x_{n,n})\big).
}
holds.}

\begin{proof} First we prove that \eq{IEG} is satisfied.

Let $x_{i,j}\in I$ for $i,j\in\{1,\dots,n\}$ and define, for $i\in\{0,\dots,n-1\}$,
\Eq{*}{
  y_i:=A^{\langle i\rangle}_{f_1,\dots,f_n}(x_{1,i+1},\dots,x_{n,i+1})
  \qquad\mbox{and}\qquad
  z_i:=A^{\langle i\rangle}_{f_1,\dots,f_n}(x_{i+1,1},\dots,x_{i+1,n}).
}
Then
\Eq{*}{
  y_i&=A_{f_1,\dots,f_n}(x_{\sigma^i(1),i+1},\dots,x_{\sigma^i(n),i+1})
   =(f_{1}+\cdots+f_n)^{-1}\bigg(\sum_{j=1}^n f_j(x_{\sigma^i(j),i+1})\bigg), \\
  z_i&=A_{f_1,\dots,f_n}(x_{i+1,\sigma^i(1)},\dots,x_{i+1,\sigma^i(n)})
   =(f_{1}+\cdots+f_n)^{-1}\bigg(\sum_{j=1}^n f_j(x_{i+1,\sigma^i(j)})\bigg).
}
Therefore, for the left and right hand sides of \eq{IEG}, we obtain the following expressions:
\Eq{*}{
  A_{f_1+\dots+f_n}(y_0,\dots,y_{n-1})
    &=(f_{1}+\cdots+f_n)^{-1}
    \bigg(\frac1n\sum_{i=0}^{n-1}\sum_{j=1}^nf_j(x_{\sigma^i(j),i+1})\bigg)\\
    &=(f_{1}+\cdots+f_n)^{-1}
    \bigg(\frac1n\sum_{\alpha=1}^{n}\sum_{\beta=1}^nf_{\sigma^{1-\beta}(\alpha)}(x_{\alpha,\beta})\bigg),\\    
  A_{f_1+\dots+f_n}(z_0,\dots,z_{n-1})&=(f_{1}+\cdots+f_n)^{-1}
    \bigg(\frac1n\sum_{i=0}^{n-1}\sum_{j=1}^nf_j(x_{i+1,\sigma^i(j)})\bigg)\\
    &=(f_{1}+\cdots+f_n)^{-1}
    \bigg(\frac1n\sum_{\alpha=1}^{n}\sum_{\beta=1}^nf_{\sigma^{1-\alpha}(\beta)}(x_{\alpha,\beta})\bigg).    
}
Thus, in order that \eq{IEG} be satisfied, it is sufficient to show that 
\Eq{Si}{
  \sigma^{1-\beta}(\alpha)=\sigma^{1-\alpha}(\beta).
}  
However, by \lem{P}, we have that $\sigma^\alpha(\alpha)=n=\sigma^\beta(\beta)$ holds for all
$\alpha,\beta\in\{1,\dots,n\}$. Applying the map $\sigma^{1-\alpha-\beta}$ to the sides of this
equation, we get
\Eq{*}{
  \sigma^{1-\beta}(\alpha)=\sigma^{1-\alpha-\beta}(n)=\sigma^{1-\alpha}(\beta),
}
which proves \eq{Si} and hence identity \eq{IEG} is also verified. 

To prove that the Gauss composition of the means 
$A^{\langle 0\rangle}_{f_1,\dots,f_n},\dots,A^{\langle n-1\rangle}_{f_1,\dots,f_n}$ is the 
quasi-arithmetic mean $A_{f_1+\dots+f_n}$, substitute $x_{i,j}:=y_i$ into \eq{IEG} where 
$y_1,\dots,y_n\in I$. Then \eq{IEG} simplifies to
\Eq{*}{
  A_{f_1+\dots+f_n}\big(A^{\langle 0\rangle}_{f_1,\dots,f_n}(y_1,\dots,y_n),
       \dots,A^{\langle n-1\rangle}_{f_1,\dots,f_n}(y_1,\dots,y_n)\big)
  =A_{f_1+\dots+f_n}(y_1,\dots,y_n).
}
Therefore $K=A_{f_1+\dots+f_n}$ is the solution of the invariance equation
\Eq{*}{
  K\circ\big(A^{\langle 0\rangle}_{f_1,\dots,f_n},\dots,A^{\langle n-1\rangle}_{f_1,\dots,f_n}\big)
  =K,
}
hence, by the unique solvability of invariance equations, $K$ is the Gauss composition of the means 
$A^{\langle 0\rangle}_{f_1,\dots,f_n},\dots,A^{\langle n-1\rangle}_{f_1,\dots,f_n}$.
\end{proof}

The following result is our main characterization theorem.

\Thm{M2}{A function $M:I^n\to I$ is a generalized quasi-arithmetic mean if and only if 
\begin{enumerate}[(i)]
 \item $M$ is a $CM$-function on $I^n$;
 \item $M$ is \emph{reflexive};
 \item $\Gamma(M^{\langle 0\rangle},\dots,M^{\langle n-1 \rangle})$ is a $CM$-function and $M$ is 
\emph{generalized bisymmetric}, that is, for all $x_{i,j}\in I$  ($i,j\in\{1,\dots,n\}$), we have
\Eq{GBS}{
  \Gamma&(M^{\langle 0\rangle},\dots,M^{\langle n-1 \rangle})\big(M^{\langle 0\rangle}(x_{1,1},
     \dots,x_{1,n}),\dots,M^{\langle n-1\rangle}(x_{n,1},\dots,x_{n,n})\big)\\
   &=\Gamma(M^{\langle 0\rangle},\dots,M^{\langle n-1\rangle})\big(M^{\langle 0\rangle}(x_{1,1},
     \dots,x_{n,1}),\dots,M^{\langle n-1 \rangle}(x_{1,n},\dots,x_{n,n})\big).
}
\end{enumerate}}

\begin{proof} If $M$ is a generalized quasi-arithmetic mean of the form $A_{f_1,\dots,f_n}$, then 
$M$ is a reflexive $CM$-function and, by \thm{M1}, the Gauss composition of the means 
$M^{\langle 0\rangle},\dots,M^{\langle n-1 \rangle}$ is the quasi-arithmetic mean $A_{f_1+\dots+f_n}$
furthermore \eq{IEG} is satisfied, which is now equivalent to \eq{GBS}.

Now assume that $M$ is a reflexive $CM$-function which satisfies \eq{GBS}. Then, using \thm{MG}, it follows
that $M$ is a $CM$-quasi-sum. Due to its reflexivity, by \lem{GQ}, we obtain that $M$ is a generalized 
quasi-arithmetic mean. 
\end{proof}


\end{document}